\documentclass{article}
\usepackage{amscd,amsfonts,euscript,amsmath,amsxtra,amssymb,amsthm}
\usepackage{mathrsfs}
\newtheorem{theorem}{Theorem}[section]

\newtheorem{proposition}[theorem]{Proposition}

\theoremstyle{definition}
\newtheorem{definition}[theorem]{Definition}
\newtheorem{example}[theorem]{Example}

\theoremstyle{remark}
\newtheorem{remark}[theorem]{Remark}

\newtheorem{convention}[theorem]{Convention}
\newtheorem{particular cases}[theorem]{Particular Cases}

\numberwithin{equation}{section}

\usepackage[all]{xy}
\def\A{{{\mathbb A}}}

\def\E{{{\mathbb E }}}

\def\Z{{{\mathbb Z }}}
\def\P{{{\mathbb P }}}

\def\AA{{{\mathscr A}}}

\def\FF{{{\mathscr F}}}

\def\II{{{\mathscr I}}}

\def\OO{{{\mathscr O}}}

\def\RR{{{\mathscr R}}}

\def\EExt{{{\mathscr E}\!xt}}

\def\FFitt{{{\mathscr F}\!itt}}

\def\Spec{{{\rm Spec \,}}}

\def\Proj{{{\rm Proj \,}}}

\def\Quot{{{\rm Quot \,}}}

\def\coker{{{\rm coker \,}}}

\def\rk{{{\rm rk \,}}}
\def\im{{{\rm im \,}}}

\def\comp{{{\mbox{\scriptsize $\circ $}}}}

\begin{document}\large 

\title{UDC 512.722+512.723;
MSC 14J60, 14D20, 14D06\;\;\;\;\;\;\;\;\;\;\;\;\;\;\;
\;\;\;\;\;\;\;\;\;\;\;\;\;\;\;\;\;\;\;\;\;\;\;\;
\;\;\;\;\;\;\;\;\;\;\;\;\;\;\;\;\;\;\;\;\;\;\;\;\;\;\;\;\;
\;\;\;\;\;\;\;\;\;\;\;\;\;\;\;\;\;\;\;\;\;\;\;\;\;\;\;\;
\;\;\;\;\;\;\;\;\;\;\;\;\;\;\;\;\;\;\;\;\;\;\;\;\;\;\\
On degeneration of surface in Fitting compactification of moduli
of stable vector bundles }
\author{N.V. Timofeeva}

\date{}

\maketitle
\begin{abstract} { The new compactification of
moduli scheme of Gieseker-stable vector bundles with the given
Hilbert polynomial on a smooth projective polarized surface
$(S,\textsf H)$, over the field $k=\bar k$ of zero characteristic,
is constructed in previous papers of the author. Families of
locally free sheaves on the surface $S$ are completed by the
locally free sheaves on the schemes which are certain
modifications of $S$. We describe the
class of modified surfaces to appear in the construction.\\
{\it Keywords:} moduli space, semistable coherent sheaves, blowup
algebra, algebraic surface. }
\end{abstract}

\section*{\rm INTRODUCTION, NOTATION AND CONVENTION} Let $S$ be a nonsingular
irreducible projective algebraic surface over an algebraically
closed field $k$ of characteristic zero. Fix an ample divisor
class ${\textsf H} \in Pic\, S.$  The symbol
 $\chi(\cdot)$ will denote, as usually, the Euler-Poincar\'{e} characteristic.
We work with the notion of (semi)stability of a coherent
torsion-free sheaf $E$ on a surface $S$ due to D.~Gieseker
\cite{Gies}.
\begin{definition} {\rm Coherent torsion-free $\OO_S$-sheaf
$E$ is {\it stable } (respectively, {\it semistable}), if for any
proper subsheaf  $F\subset E$ for } $m\gg 0$
$$
\frac{\chi(F(m{\textsf H}))}{\rk (F)}<\frac{\chi(E(m{\textsf
H}))}{\rk (E)},\;\; (\mbox{\rm respectively,} \;\;
\frac{\chi(F(m{\textsf H}))}{\rk (F)}\leq\frac{\chi(E(m{\textsf
H}))}{\rk (E)}\;).
$$
\end{definition}

It is well-known \cite{EG} that the structure of moduli space for
semistable sheaves depends strongly on the choice of polarization
$\textsf H$. The Gieseker and Maruyama moduli scheme for
semistable torsion-free sheaves on the surface $S$, with  Hilbert
polynomial $P(m)=\chi(E(m{\textsf H}))$ with respect to the class
$\textsf H$, is denoted by $\overline M.$ It is known \cite{Mar}
that it is a projective scheme of finite type over $k$. The points
corresponding to locally free sheaves (vector bundles) constitute
Zariski-open subscheme $M_0$ of $\overline M$. Assume that
$\overline M$ is a fine moduli space. Then there is a trivial
family of surfaces $\Sigma :=\overline M\times S\longrightarrow
\overline M$ carrying the universal family of stable sheaves $\E$.
In \cite{Tim0, Tim1} a projective scheme $\widetilde M$ and a
nontrivial flat family of (possibly reducible) schemes $\widetilde
\Sigma\stackrel{\widetilde \pi}{\longrightarrow}\widetilde M$
endowed with the family of locally free sheaves $\widetilde \E$,
are constructed. In \cite{Tim2} the analogous construction (flat
families of schemes $\widetilde \Sigma_i \stackrel{\widetilde
\pi_i}{\longrightarrow}\widetilde B_i$ over \'{e}tale
neighborhoods  $\widetilde B_i
\stackrel{\acute{e}tale}{\longrightarrow }\widetilde M$, endowed
with locally free sheaves $\widetilde \E_i$) is done for the
coarse case. The scheme $\widetilde M$ contains Zariski-open
subscheme $\widetilde M_0$ isomorphic to $M_0$. Moreover, in both
(fine and coarse) cases there is a birational morphism of
compactifications $\varphi: \widetilde M \to \overline M$ such
that $\varphi|_{\widetilde M_0}$ is isomorphism.

\begin{convention} {\rm We work in the notations for the fine case.
Although, all the considerations will be valid for the coarse case
as well.}
\end{convention}

The birational morphism $\varphi: \widetilde M \to \overline M$
constructed in papers  \cite{Tim0, Tim1},\linebreak establishes
the correspondence among pairs $(\widetilde S, \widetilde E)\in
\widetilde M$ and $\varphi(\widetilde S, \widetilde E)=(S,E)\in
\overline M.$ Let $\widetilde y \in \widetilde M$ and $\varphi
(\widetilde y)=y$. Hence  we mean that the fibre $\pi^{-1}(y)=S$
of family $\Sigma$ is the image of the fibre $\widetilde
\pi^{-1}(\widetilde y)=\widetilde S$, and the coherent sheaf  $E$
on a fibre  $S$ is the image of the vector bundle $\widetilde E$
on the fibre $\widetilde S$.

Recall the follow\-ing definition.
\begin{definition} {\rm (sheaf analogue of given in \cite[20.2]{Eisenbud}) Let $X$ be
an algebraic scheme, $F_0,\; F_1$ locally free $\OO_X$ -- sheaves,
$\psi: F_1\to F_0$ -- $\OO_X$-module homomorphism. Denote $\FF=
\coker\, \psi, $ $r_0=\rk \,F_0.$ {\it The sheaf of zeroth Fitting
ideals of $\OO_X$- module $\FF$} is defined as
$$\FFitt^0(\FF)=\im (\psi':\bigwedge^{r_0} F_1 \otimes
\bigwedge^{r_0} F_0^{\vee}\to \OO_X),$$ where $\psi'$ is the
associate morphism for $\psi.$}\end{definition}

The aim of the present paper is to investigate the structure of
fibres of the morphism $\widetilde \pi: \widetilde \Sigma \to
\widetilde M$ in general case. In \cite{Tim1} it is proven that
the  fibre at a general point $\widetilde y \in \widetilde M_0$ is
isomorphic to $S$ and one component of the  fibre at a special
point $\widetilde y \in \widetilde M \backslash \widetilde M_0$ is
isomorphic to the blowup of $S$ in the sheaf of ideals $\FFitt^0
\EExt^1(E, \OO_S).$ Now we give a description for the whole of the
scheme $\widetilde S.$

As proven in \cite[Proposition 3.1]{Tim1}, the scheme $\widetilde
\Sigma $ is given by the blowup $\sigma \!\!\! \sigma: \widetilde
\Sigma \to \widetilde M \times S$ of the trivial family
$\widetilde M\times S $ in the sheaf of ideals $\FFitt^0
\EExt^1_{\Delta}(\OO_{\widetilde M}\boxtimes \E|_{\Delta},
\OO_{\Delta})$ with $\Delta$ being the scheme-theoretic closure of
the image of diagonal immersion $\widetilde M_0 \times S
\hookrightarrow \widetilde M_0 \times M_0 \times S$ in the product
$\widetilde M \times \Sigma$. The scheme $\Delta$ is isomorphic to
$\widetilde M \times S$. This means that to describe fibres of the
projection $\widetilde \pi$ it is enough to investigate fibres of
the composite morphism $\pi \comp \sigma\!\!\!\sigma: \widetilde
\Sigma \stackrel{\sigma \!\!\! \sigma }{\longrightarrow}
\widetilde M \times S \stackrel{\pi}{\longrightarrow} \widetilde
M.$

The paper is organized as follows. In the first section we give
some observations about  the structure of singularities of
semistable torsion-free coherent sheaves on a nonsingular surface.
In the second section we construct a flat 1-dimensional family of
coherent sheaves such that it contains $E$ and its general sheaf
is locally free. At last, in the third section we derive the
scheme structure of $\widetilde S$.
\begin{convention} {\rm As usually we assume that $\overline M$ is irreducible. If
not, consider each of its irreducible components containing
locally free sheaves. We restrict ourselves to the case when
$\overline M$ contains at least one locally free sheaf.}
\end{convention}

\begin{convention} {\rm In the whole of the text we omit the subscripts in $\EExt$'s whenever no
ambiguity occur: for example, we replace $\EExt^1_{S}(E,\OO_S)$ by
$\EExt^1(E,\OO_S).$ When working with Artinian sheaves the length
$l(\varkappa)$ of a sheaf $\varkappa$ denotes its
Euler-Poincar\'{e} characteristic: $l(\varkappa )= \chi
(\varkappa)$. The same is for zero dimensional subscheme $Z\subset
S$:\linebreak $l(Z)=l(\OO_Z)=\chi(\OO_Z).$  As usually, we denote
the Grothendieck scheme of length $l$ zero dimensional quotients
of $\OO_S$-sheaf $F$ on $S$ as $\Quot^l F.$ The point in the
Grothendieck scheme, corresponding to the quotient $q:
F\twoheadrightarrow \varkappa$, is denoted as  $q$.}
\end{convention}

Let $r$ be the rank of coherent sheaves $E$ with Hilbert
polynomial $P(m)$. The final result of this paper is given by the
following theorem.
\begin{theorem}{The  fibre of the family
$\widetilde \pi: \widetilde \Sigma \to \widetilde M$ at the point
$\widetilde y\in \widetilde M$

i) is isomorphic to $S$ if $\widetilde y \in \widetilde M_0$, or

ii) contains in the class of all $\Proj \bigoplus_{s\ge 0}
(\II_S[t]+(t))^s/(t^{s+1})$, with $\II_S=\linebreak \FFitt^0
\EExt^2(\varkappa, \OO_S)$ where $\varkappa$ denotes the length
$l$ Artinian sheaf which is a quotient of the direct sum
$\OO_S^{\oplus r}$, $l\le c_2$, if $\widetilde y\in \widetilde M
\backslash \widetilde M_0$. }
\end{theorem}

{\bf Acknowledgements.} The author expresses her deep and sincere
gratitude to M.~Reid in  Mathematics Research Centre of the
University of Warwick, Great Britain, V.\,V.~Shokurov and D.\,B.~
Kaledin in V.\,A.~Steklov Institute of RAS, Moscow, Russia, and M.
E. Sorokina in Yaroslavl State Pedagogical Unversity, Yaroslavl,
 Russia, for discussions and comments.

\section{\rm SINGULARITY OF  $E$}

In \cite{Tim1} we proved that the main (dominating $S$) component
$\widetilde S_0$ of $\widetilde S$ is the blowup of $S$ in the
sheaf of ideals $\FFitt^0 \EExt^1(E,\OO_S)$. In this section we
investigate the class of sheaves who appear as $\EExt^1(E,\OO_S)$
for various semistable $E$.

For any torsion-free sheaf $E$ there is exact triple
\begin{equation}\label{mex}
0\to E \to E^{\vee \vee} \to \varkappa \to 0.
\end{equation}
Since $S$ is nonsingular surface, the double dual sheaf $E^{\vee
\vee}$ is locally free and the cokernel $\varkappa$ is Artinian
sheaf.
\begin{definition} {\rm In (\ref{mex}) the cokernel $\varkappa$ is said to be
{\it a singularity sheaf } of $E$. When necessary, we reflect this
fact in the notation } $\varkappa_E:=E^{\vee \vee}/E.$
\end{definition}

The form of Hilbert polynomial of $E$ is determined by the
geometry of the surface $S$, choice of polarization $\textsf H$,
rank $r$ and Chern classes $c_1, c_2$ of the sheaf $E$. In any
case, all possible $l$'s are bounded by the inequality $0\le l \le
c_2$. But for $c_2$ fixed, and various $S, {\textsf H}, r, c_1$,
the collections of possible $l$'s can be different.  However,
there is no explicit description of such $l$'s up to now.

Now recall the notion of slope-(semi)stability ascending to
D.~Mumford and F.~Takemoto. We use the following definition.

\begin{definition} {\rm Coherent torsion-free sheaf $E$ is
{\it slope-(semi)stable} with respect to the polarization $\textsf
H$ if for any proper subsheaf $F \subset E$ the following holds:}
\begin{equation*}
\frac{c_1(F)\cdot {\textsf H}}{\rk(F)}< \frac{c_1(E)\cdot {\textsf
H}}{\rk(E)}, \quad (\mbox{\rm respectively, } \frac{c_1(F)\cdot
{\textsf H}}{\rk(F)}\le \frac{c_1(E)\cdot {\textsf H}}{\rk(E)}).
\end{equation*}
\end{definition}
Semistability implies slope-semistability, slope-stability implies
stability.

The following simple remark shows that the possible $l$'s may
cover the interval $[0, c_2]\subset \Z$ not completely.

\begin{remark} {\rm For $r,c_1,{\textsf H}$ such that $r$ and
$(c_1 \cdot {\textsf H})$ coprime, the slope-semi\-stability
implies slope-stability. In this case all semistable sheaves are
slope-stable. Let the  invariants $r, c_1, c_2, l,$ $l\le c_2$, be
such that there are no slope-stable sheaves $E$ with $\rk (E)=r,
c_1(E)=c_1, c_2(E)=c_2$ but there exist at least one slope-stable
vector bundle $F$ of the same rank and first Chern class but
$c_2(F)=c_2-l$. Then for any Artinian sheaf $\varkappa$ of length
$l$ the slope-stable kernel $E$ of the morphism
$F\twoheadrightarrow \varkappa$ must not exist. This means that
length-$l$-sheaves do not appear as cokernels of exact triples
(\ref{mex}).}
\end{remark}

\begin{example} {\rm These effects can be observed even in
the classical case $S=\P^2$, $r=2,$ $c_1=-1,$ $c_2=2$. The
corresponding moduli variety of semistable (=stable=slope-stable)
coherent sheaves contains nonempty locus $M_0$ of locally free
sheaves \cite[Ch. II, 4]{OSS}. Assume that there is a stable sheaf
$E$ with $l=2$. Then for $E^{\vee \vee}$ compute Bogomolov's
discriminant: $\Delta (E^{\vee \vee})=2rc_2(E^{\vee
\vee})-(r-1)c_1^2(E^{\vee \vee})=-1$. This contradicts the
slope-stability of $E^{\vee \vee}$ and shows that for nonlocally
free $E$ there is only one possibility $E^{\vee \vee}/E\cong k.$}
\end{example}

\begin{definition} {\rm The Artinian sheaf $\varkappa$ is said to be
$(S, {\textsf H}, r, P(m))$-{\it ad\-miss\-ible} if there is an
exact $\OO_S$-triple (1.1) where coherent sheaf $E$ of rank $r$
with Hilbert polynomial $P(m)$ is semistable with respect to the
polarization $\textsf H$. }
\end{definition}

Applying functor $\EExt^{\bullet}(-, \OO_S)$ to (\ref{mex}) one
gets immediately
\begin{equation*}
\EExt^1(E, \OO_S)=\EExt^2(\varkappa, \OO_S)
\end{equation*}
and we have

 \begin{proposition}\label{classext} {Class of all sheaves of
 ideals
 $\FFitt^0\EExt^1(E,\OO_S)$ for all semistable $E$ of fixed
rank $r$ and Hilbert polynomial $P(m)$, contains in  the class of
sheaves as follows
\begin{equation*} \FFitt^0
\EExt^2(\varkappa, \OO_S)
\end{equation*}
for all $q:\OO_S^{\oplus r}\twoheadrightarrow \varkappa$, $q\in
\Quot^l \OO_S^{\oplus r},$ $l\le c_2.$ }
\end{proposition}

\section{\rm REMOVABILITY OF SINGULARITY}

We prove (with respect to convention 0.3) that any semistable
coherent sheaf can be include into some 1-dimensional flat family
of sheaves  such that the general sheaf of the family is locally
free.

Take $m\gg 0$ such that $E(m{\textsf H})$ is globally generated.
Consider the vector space $V\cong H^0(S, E(m{\textsf H}))$ and
Grothendieck's Quot-scheme  $Quot=\Quot^{P(m)} (V\otimes
\OO_S(-m{\textsf H}))$ parameterizing quotient sheaves
\begin{equation} V\otimes \OO_S(-m{\textsf H})\twoheadrightarrow E
\end{equation} with Hilbert polynomial equal to
$P(m)=\chi(E(m{\textsf H}))$. We work as usually with the
quasiprojective subscheme $Q\subset Quot$ constituted by all
quotients (2.1) with $E$ semistable and with isomorphism $H^0(S,
E(m{\textsf H}))\cong V.$ Grothendieck's scheme $Quot$ carries the
universal quotient sheaf $\E_{Quot}$, let $\E_Q$ be its
restriction onto $Q$. Let $Q_0$ be the open subscheme of $Q$ whose
points correspond to locally free quotient sheaves $E$.

The further consideration contains the usage of Bertini's theorem
and one needs the smoothness of $Q.$ If it is not so, replace it
by any  smooth resolution, for example due to H. Hironaka
\cite{Hiron}.

Since $Q$ is quasiprojective, there is a projective space $\P^N$
together with locally closed immersion $i: Q \hookrightarrow
\P^N.$ Fix any point $q$ in the image $i(Q)$  such that $q$
corresponds to $E$. Consider the set of all hyperplanes in $\P^N$
passing through $q$. Choose hyperplane $H_1$ such that the
intersection scheme $Q_{(1)}=i(Q) \cap H_1$ containing $q$ and
meeting $i(Q_0)$ is irreducible and nonsingular. It is possible by
Bertini's theorem \cite[Ch. III, Corollary 7.9]{Hart1}. Clearly,
$i(Q_0)\cap Q_{(1)}$ forms an open subset in $Q_{(1)}.$ Now repeat
the procedure replacing $Q$ by $Q_{(1)}.$ Iterating the process
one comes to $Q_{(d)}$ being irreducible curve for some $d>0$.

\section{\rm STRUCTURE OF MODIFIED SURFACES}
Here we derive the scheme structure of surfaces $\widetilde S$ as
projective spectra of appropriate algebras. As usually $\OO_S$
denotes the structure sheaf of the surface $S$ and let $\II_S
\subset \OO_S$ be the sheaf of ideals to be blown up in $S$.

As direct computation with blowup equations shows, the scheme
$\widetilde S$ can carry quite sophisticated structure. The main
component $\widetilde S_0$ admits singularities and each of other
components can carry nonreduced scheme structure.

\begin{example} {\rm By locality one can replace the original nonsingular
surface by the affine subset $U\cong \A^2=\Spec [x,y]$. Take
$I_S=\Gamma (U,\II_S)=(x^2,y)$. Consider the trivial family
$T_U=U\times \Spec k[t]=\Spec k[x,y,t]$ with natural projection
$\pi:T_U \to \Spec k[t].$ We blow up $T_U$ in the nonreduced point
with ideal  $I_S$ in the fibre over $b_0=\{t=0\}$. This is
equivalent to the blowup $\sigma \!\!\! \sigma: \widetilde T_U \to
T_U$ of the point with ideal $I=(x^2,y,t)$ on $T_U$. The  scheme
$\widetilde T_U$ is given in the direct product $T_U \times
\P^2\cong \A^3 \times \P^2$ for $\P^2=\Proj k[u,v,w],$ by usual
blowup relations: $x^2:y:t=u:v:w$. The exceptional divisor
$\textsf{E}$ of the blowup morphism $\sigma\!\!\!\sigma:
\widetilde T_U \to T_U$ carries a nonreduced ("double") structure;
its equations are $x^2=y=t=0$. The fibre $\widetilde S$ of
composite morphism $\pi\comp \sigma\!\!\!\sigma$ over $b_0$
consists of two components: $\widetilde S=\widetilde S_0 \cup
\widetilde S_1$. In our case $\widetilde S_0$ has a quadratic
singularity; in the affine chart $v \ne 0,$ $z:=u/v$ its equation
in the neighborhood of singular point is $x^2=yz$. The component
$\widetilde S_1=\textsf{E}$.}
\end{example}

For the general consideration form a polynomial extension
$\OO_S[t]$ for $t$ trans\-cendental over $\OO_S$, let $(t)\subset
\OO_S[t]$ be principal ideal sheaf. Set $\II:=\II_S[t]+(t)\subset
\OO_S[t]$. Set as well $T:=\Proj \OO_S[t]$, and $\OO_Z:=\OO_T
/\II$. Clearly, $T=\Spec k[t] \times S$. Denote $B:=\Spec k[t]$,
and the zero point on the base $B$ is $b_0=\{t=0\}$. Let $\pi:
T\to B$ be the projection induced by the $\OO_S[t]$ - algebra
morphism $k[t] \to \OO_S[t]$. The latter morphism is obtained by
the extension of the structure morphism $k \to \OO_S.$

Form a graded sheaf algebra $\AA:= \bigoplus_{s\ge 0} \II^s$, and
$\widetilde T:=\Proj \AA$. There is a projective (blowup) morphism
$\sigma\!\!\!\sigma : \widetilde T \to T$ induced by the natural
$\OO_S[t]$-algebra morphism $\OO_S[t] \to \AA$ onto the zero
graded component.

\begin{proposition}\label{fibre} { The fibre of composite morphism $\sigma\!\!\!\sigma \comp \pi$ at
$b_0 \in B$ equals
$$
\widetilde S=(\Proj \AA)\times _T \pi^{-1} (b_0)=\Proj \bigoplus
_{s\ge 0}\II^s / (t^{s+1}).$$}
\end{proposition}
\begin{example}
{\rm Take $\II_S={\mathfrak m}_{p}$ -- the sheaf of maximal ideals
corresponding to a reduced point $p\in S=\pi^{-1}(b_0)$. After the
restriction to appropriate affine neighborhood $U$ one has
$I=(x,y,t),$ henceforth $A=\Gamma(U,\AA)=\bigoplus_{s\ge
0}(x,y,t)^s,$ $T=\Spec k[x,y,t],$ $\widetilde T= \Proj
\bigoplus_{s \ge 0} (x,y,t)^s$. For the ex\-ceptional divisor
$\textsf{E}=\sigma\!\!\!\sigma ^{-1}(p,b_0)$ one has
$$
\textsf{E}=(\Proj \bigoplus_{s\ge 0}(x,y,t)^s) \times _T
(p,b_0)=\Proj \bigoplus_{s\ge 0} (x,y,t)^s/(x,y,t)^{s+1}.
$$

The fibre $\widetilde S=(\sigma\!\!\!\sigma \comp \pi)^{-1} (b_0)
$  is given by the fibered product
\begin{equation}\label{square}\xymatrix{\widetilde T \ar@{>}[r] ^{\sigma\!\!\!\sigma} &T\\
\widetilde S \ar@{^{(}->}[u] \ar[r]^{\sigma } & S \ar@{^{(}->}[u]}
\end{equation}
with $\sigma $ be the restriction of $\sigma\!\!\!\sigma $ onto
the fibre $\widetilde S$ over $b_0$. Passing to algebras one has
that the algebra $R=\bigoplus_{s\ge 0} R_s$ for $\widetilde S$ is
the coproduct  of graded algebras given by
\begin{equation*}
R_s=(x,y,t)^s \coprod _{ k[x,y,t]} k[x,y].
\end{equation*}
The push-out diagram of $k[x,y,t]$-modules
\begin{equation}\label{pushout}\xymatrix{&0 \ar[d]& 0\ar[d]\\
&(t) \ar[d] \ar@{=}[r]& (t)\ar[d]^{\cdot t^s}\\
0\ar[r]& k[x,y,t] \ar[r]^{\cdot t^s} \ar[d]& (x,y,t)^s \ar[r]
\ar[d]& C \ar[r] \ar@{=}[d]& 0\\
0\ar[r]& k[x,y] \ar[r] \ar[d] & R_s \ar[r] \ar[d]& C \ar[r]&0\\
& 0& 0\\ }
\end{equation}
gives the explicit form of $R_s$:
\begin{equation}\label{rn} R_s=(x,y,t)^s/(t^{s+1})\end{equation}

The universal property of $R$ as a coproduct is checked
immediately.

The inclusion of the exceptional divisor $\textsf{E}$ into
$\widetilde S$ is defined by the epi\-morphism of algebras
\begin{equation} \label{onto}\bigoplus_{s\ge 0} (x,y,t)^s
/(t^{s+1})\twoheadrightarrow \bigoplus_{s\ge 0} (x,y,t)^s
/(x,y,t)^{s+1}.\end{equation} As well the inclusion of the main
component $\widetilde S_0 \hookrightarrow \widetilde S$ is defined
by the epimorphism}
\begin{equation}\label{onto1}
\bigoplus_{s\ge 0}(x,y,t)^s /(t^{s+1})\twoheadrightarrow
\bigoplus_{s\ge 0} (x,y)^s.
\end{equation}
\end{example}
\begin{proof}[Proof of proposition \ref{fibre}] Let $Z\subset T$ be
zero dimensional subscheme defined by the sheaf of ideals $\II$.
The exceptional divisor $\textsf{E}$ of the blowup
$\sigma\!\!\!\sigma$ is given by
$$\textsf{E}:=(\Proj \AA)\times _T Z=\Proj \bigoplus_{s\ge 0} \II^s /
\II^{s+1}.$$ The fibered square (\ref{square}) relates to the
(sheaf) algebra
\begin{equation}\label{rgrad}\RR=\bigoplus_{s\ge 0}\RR_s\end{equation} as a coproduct
\begin{equation*} \RR_s= \II^s
\coprod_{\OO_S[t]} \OO_S,\;\; s\ge 0.
\end{equation*}
For (\ref{pushout}) and (\ref{rn}) one has straightforward
generalizations and
\begin{equation}\label{rrn} \RR_s=\II^s/(t^{s+1}).
\end{equation}
Analogously to (\ref{onto}), (\ref{onto1}) there are the
inclusions of exceptional divisor $\textsf{E}$ and of the main
component $\widetilde S_0$ defined by the sheaf algebra
epimorphisms
\begin{eqnarray}
\bigoplus_{s\ge 0} \II^s /(t^{s+1})&\twoheadrightarrow&
\bigoplus_{s\ge 0} \II^s /\II^{s+1}, \nonumber\\ \bigoplus_{s\ge
0}\II^s /(t^{s+1})&\twoheadrightarrow& \bigoplus_{s\ge 0} \II_S^s
\nonumber
\end{eqnarray}
respectively. Hence, $\widetilde S=\Proj \RR$ for $\RR$ given by
(\ref{rgrad}), (\ref{rrn}). \end{proof}
\begin{remark} {\rm Let $\textsf{E}_0$ be the exceptional divisor of the main
component $\widetilde S_0$, i.e. exceptional divisor for the
blowup morphism $\sigma_0=\sigma|_{\widetilde S_0}: \widetilde S_0
\to S$. It is easily seen that the main component $\widetilde
S_0=\Proj \bigoplus_{s\ge 0}\II_S^s$ of special fibre $\widetilde
S$ meets $\textsf{E}$ precisely at $\textsf{E}_0=\Proj
\bigoplus_{s\ge 0}\II_S^s/\II_S^{s+1}.$}
\end{remark}


\vspace{15mm}
\begin{flushleft}
Nadezda V. Timofeeva\\
Department of Algebra and Mathematical Logic, \\
Yaroslavl' State University,\\
Sovetskaya str., 14, Yaroslavl', 150000\\
Russia\\
{\it E-mail:} ntimofeeva@list.ru
\end{flushleft}

\end{document}